\documentclass{elsart}
\usepackage{amsfonts}
\usepackage{amsmath}

\begin{document}

\begin{frontmatter}%

\title{On Better Approximation Order for the Nonlinear Baskakov Operator
	of Maximum Product Kind}%

\author{Sezin \c{C}\.{I}T\corauthref{cor}}
\corauth[cor]{Corresponding author}%

\address{Department of Mathematics, Faculty of Science, Gazi University,
	Ankara, Turkey}%

\ead{sezincit@gazi.edu.tr}%

\author{Og\"{u}n DO\u{G}RU}%

\address{Department of Mathematics, Faculty of Science, Gazi University,
	Ankara, Turkey}%

\ead{ogun.dogru@gazi.edu.tr}%

\begin{abstract}
Using maximum instead of sum, nonlinear Baskakov operator of maximum product
kind is introduced by Bede et al. \cite{max-prod Baskakov}. The present
paper deals with the approximation processes for this operator. Especially
in \cite{max-prod Baskakov}, it was indicated that the order of
approximation of this operator to the function $f$ under the modulus is $%
\frac{\sqrt{x\left( 1+x\right) }}{\sqrt{n}}$ and it could not be improved
except for some subclasses of functions. Contrary to this claim under some
circumstances, we will show that a better order of approximation can be
obtained with the help of classical and weighted modulus of continuities. 
\end{abstract}%


\end{frontmatter}%

\section{Introduction}

For $f\in C[0,\infty ),$ the classical Baskakov operators defined as%
\begin{equation*}
	V_{n}(f;x)=\frac{1}{(1+x)^{n}}\sum_{k=0}^{\infty }\binom{n+k-1}{k}\left( 
	\frac{x}{1+x}\right) ^{k}f\left( \frac{k}{n}\right) ,
\end{equation*}%
where $x\in \left[ 0,\infty \right) $, $n\in \mathbb{N}$, were introduced in 
\cite{Baskakov}.

The construction logic of nonlinear maximum product type operators using the
maximum instead of the sum is based on the studies \cite{max-prod-shepard}, 
\cite{approx by pseudo} and \cite{Gal} (for details, see also, \cite%
{max-prod kitap}).

There are some other remarkable articles like \cite{max-prod-berns}, \cite%
{max-prod-Szasz}, \cite{max-prod-berns-szazs}, \cite{Max-prod-
	Bleimann-Butz-Hahn} and \cite{max-prod Baskakov} that we will remind you in
chronological order that various maximum product type nonlinear operators
was introduced and their approximation and rate of convergence properties
were investigated. Also, some statistical approximation properties of
maximum product type operators was given by Duman in \cite{Duman}.

Especially, in \cite{max-prod Baskakov}, the approximation properties, rate
of convergence and shape preserving properties of Baskakov operator of
maximum product kind are examinated.

At this point let us recall the following well known concept of classical
modulus of continuity:%
\begin{equation}
	\omega \left( f,\delta \right) =\max \left\{ \left\vert f\left( x\right)
	-f\left( y\right) \right\vert ;\text{ }x,y\in I,\text{ }\left\vert
	x-y\right\vert \leq \delta \right\} .  \label{1.1}
\end{equation}%
The order of approximation for the maximum product type Baskakov operator
can be found in \cite{max-prod Baskakov} by means of the modulus of
continuity as \\ $\omega \left( f;\sqrt{x\left( 1+x\right) /n}\right) $. Also,
Bede et al. indicated that the order of approximation under the modulus was $%
\sqrt{x\left( 1+x\right) /n}$ and it could not be improved except for some
subclasses of functions (see, for details, \cite{max-prod Baskakov}).

Contrary to this claim under some circumstances, we will show that a better
order of approximation can be obtained with the help of classical and
weighted modulus of continuities.

Notice that, in \cite{better Szasz}, \cite{better BHH}, \cite{better
	bernstein} we showed that rates of approximation of some max-product type
operators can be improved.

\section{The Concept of Nonlinear Maximum Product Operators}

Before the given of the main results, we will recall basic definitions and
theorems about nonlinear operators given in \cite{max-prod-berns-szazs}, 
\cite{max-prod kitap} and \cite{approx by pseudo}.

The set of non-negative real numbers, $%
\mathbb{R}
_{+}$ is $\vee $ (maximum) and $\cdot $ (product). Let's consider the
operations. Then $\left( 
\mathbb{R}
_{+},\vee ,\cdot \right) $ has a semi-ring structure, which is called a
maximum-product algebra.

Let $I\subset \mathbb{R}$ be bounded or unbounded interval, and 
\begin{equation*}
	CB_{+}\left( I\right) =\left\{ f:I\rightarrow \mathbb{R}_{+}:f\text{
		continuous and bounded on }I\text{ }\right\} .
\end{equation*}

Let us take the general form of $L_{n}:CB_{+}(I)\rightarrow CB_{+}(I),$ as 
\begin{equation*}
	L_{n}\left( f\right) (x)=\bigvee\limits_{i=0}^{n}K_{n}\left( x,x_{i}\right)
	f(x_{i})\text{ }
\end{equation*}%
or%
\begin{equation*}
	L_{n}\left( f\right) (x)=\bigvee\limits_{i=0}^{\infty }K_{n}\left(
	x,x_{i}\right) f(x_{i}),
\end{equation*}%
where $n\in \mathbb{N}$, $f\in CB_{+}(I)$, $K_{n}\left( .,x_{i}\right) \in
CB_{+}\left( I\right) $ and $x_{i}\in I$, for all $i.$ These operators are
nonlinear, positive operators and moreover they satisfy the following
pseudo-linearity condition of the form%
\begin{equation*}
	L_{n}\left( \alpha f\vee \beta g\right) \left( x\right) =\alpha
	\,L_{n}\left( f\right) \left( x\right) \vee \beta \,L_{n}\left( g\right)
	\left( x\right) ,\forall \alpha ,\beta \in \mathbb{R}_{+},\text{ }f,g\in
	CB_{+}(I).
\end{equation*}

In this section, we present some general results on these kinds of operators
which will be used later.

\begin{lem}
	\label{Lemma2.1} \cite{max-prod-berns-szazs} Let $I\subset \mathbb{R}$ be
	bounded or unbounded interval,%
	\begin{equation*}
		CB_{+}\left( I\right) =\left\{ f:I\rightarrow \mathbb{R}_{+}:f\text{
			continuous and bounded on }I\text{ }\right\} ,
	\end{equation*}%
	and $L_{n}:CB_{+}(I)\rightarrow CB_{+}(I),$ $n\in \mathbb{N}$ be a sequence
	of operators satisfying the following properties:
	
	$(i)$ If $f,g\in CB_{+}\left( I\right) $ satisfy $f\leq g$ then $L_{n}\left(
	f\right) \leq L_{n}\left( g\right) $ for all $n\in \mathbb{N}$.
	
	$(ii)$ $L_{n}\left( f+g\right) \leq L_{n}\left( f\right) +L_{n}\left(
	g\right) $ for all $f,g\in CB_{+}\left( I\right) .$
	
	Then for all $f,g\in CB_{+}\left( I\right) ,$ $n\in \mathbb{N}$ and $x\in I$
	we have%
	\begin{equation*}
		\left\vert L_{n}\left( f\right) \left( x\right) -L_{n}\left( g\right) \left(
		x\right) \right\vert \leq L_{n}\left( \left\vert f-g\right\vert \right)
		\left( x\right) .
	\end{equation*}
\end{lem}

After this point, let us denote the monomials $e_{r}(x):=x^{r},$ $r\in 
\mathbb{N}_{0}.$ First three monomials are also called as Korovkin test
functions.

\begin{cor}
	\label{Corollary2.1} \cite{max-prod-berns-szazs} Let $L_{n}:CB_{+}(I)%
	\rightarrow CB_{+}(I),$ $n\in \mathbb{N}$ be a sequence of operators
	satisfying the conditions $(i),(ii)$ in Lemma \ref{Lemma2.1} and in addition
	being positive homogenous. Then for all $f\in CB_{+}(I)$, $n\in \mathbb{N}$
	and $x\in I$ we have 
	\begin{equation*}
		\left\vert L_{n}\left( f\right) \left( x\right)-f(x)\right\vert \leq \left[ 
		\frac{1}{\delta }L_{n}\left( \varphi _{x}\right) \left( x\right)
		+L_{n}\left( e_{0}\right) \left( x\right) \right] \omega \left( f,\delta
		\right)+f\left( x\right) \,\left\vert L_{n}\left( e_{0}\right) \left(
		x\right) -1\right\vert ,
	\end{equation*}%
	where $\omega \left( f,\delta \right) $ is the classical modulus of
	continuity defined by (\ref{1.1}), $\delta >0,$ $e_{0}\left( t\right) =1,$ $%
	\varphi _{x}\left( t\right) =\left\vert t-x\right\vert $ for all $t\in I,$ $%
	x\in I,$ and if $I$ is unbounded then we suppose thet there exists $%
	L_{n}\left( \varphi _{x}\right) \left( x\right) \in \mathbb{R}_{+}\cup
	\left\{ \infty \right\} ,$ for any $x\in I,$ $n\in \mathbb{N}$.
\end{cor}

A consequence of Corollary \ref{Corollary2.1}, we have the following:

\begin{cor}
	\label{Corollary2.2} \cite{max-prod-berns-szazs} Suppose that in addition to
	the conditions in Corollary \ref{Corollary2.1}, the sequence $\left(
	L_{n}\right) _{n}$ satisfies $L_{n}\left( e_{0}\right) =e_{0},$ for all $%
	n\in \mathbb{N}$. Then for all $f\in CB_{+}\left( I\right) $, $n\in \mathbb{N%
	}$ and $x\in I$ we have 
	\begin{equation*}
		\left\vert L_{n}\left( f\right) \left( x\right) -f(x)\right\vert \leq \left[
		1+\frac{1}{\delta }L_{n}\left( \varphi _{x}\right) \left( x\right) \right]
		\,\omega \left( f,\delta \right)
	\end{equation*}%
	where $\omega \left( f,\delta \right) $ is the classical modulus of
	continuity defined by (\ref{1.1}) and $\delta >0.$
\end{cor}

\section{Nonlinear Baskakov Operator of Maximum Product Kind}

\bigskip In the classical Baskakov operator, the sum operator $\sum $ is
replaced by the maximum operator $\bigvee $, and introduced by Bede et al.
in \cite{max-prod Baskakov}. So, nonlinear Baskakov operator of maximum
product kind is defined as 
\begin{equation}
	V_{n}^{(M)}(f)(x):=\frac{\bigvee\limits_{k=0}^{\infty }b_{n,k}\left(
		x\right) \,f\left( \frac{k}{n}\right) }{\bigvee\limits_{k=0}^{\infty
		}b_{n,k}\left( x\right) },  \label{2.1}
\end{equation}%
where $b_{n,k}\left( x\right) =\binom{n+k-1}{k}\frac{x^{k}}{\left(
	1+x\right) ^{n+k}},$ $f\in C[0,\infty ),$ $x\in \left[ 0,\infty \right) $, $%
n\in \mathbb{N}$.

\begin{rem}
	In \cite{max-prod Baskakov}, the approximation and shape preserving
	properties of the operator $V_{n}^{(M)}(f)(x)$ are also examinated.
\end{rem}

\begin{lem}
	\cite{max-prod Baskakov} $1)$ It is easy to see that the nonlinear Baskakov
	max-product operator satisfy the conditions $(i),$ $(ii)$ of Lemma \ref%
	{Lemma2.1}. In fact, instead of $(i)$ it also satisfies the following
	stronger condition:%
	\begin{equation*}
		V_{n}^{(M)}\left( f\vee g\right) \left( x\right) =\,V_{n}^{(M)}\left(
		f\right) \left( x\right) \vee \,V_{n}^{(M)}\left( g\right) \left( x\right) ,%
		\text{ }f,g\in CB_{+}(I),\text{ }I=\left[ 0,\infty \right) .
	\end{equation*}%
	Indeed, taking into consideration of the equality above, for $f\leq g,$ $%
	f,g\in CB_{+}(I),$ it easily follows $V_{n}^{(M)}\left( f\right) (x)\leq
	V_{n}^{(M)}\left( g\right) (x).$
	
	$2)$ In addition to this, it is immadiate that the nonlinear Baskakov
	max-product operator is positive homogenous, that is $V_{n}^{(M)}\left(
	\lambda f\right) =\lambda V_{n}^{(M)}\left( f\right) $ for all $\lambda \geq
	0.$
\end{lem}

\begin{lem}
	\label{Lemma3.1} \cite{max-prod Baskakov} For any arbitrary bounded function 
	$f:\left[ 0,\infty \right) \rightarrow \mathbb{R}_{+}$, max-product operator 
	$V_{n}^{(M)}(f)(x)$ is positive, bounded, continuous and satisfies $%
	V_{n}^{(M)}(f)(0)=f\left( 0\right) ,$ for all $n\in \mathbb{N}$, $n\geq 3.$
\end{lem}

\section{Auxiliary Results}

\begin{lem}
	\label{Lemma4.1} \cite{max-prod Baskakov} Let $n\in \mathbb{N},$ $n\geq 2.$
	We have 
	\begin{equation*}
		\bigvee\limits_{k=0}^{\infty }b_{n,k}\left( x\right) =b_{n,j}(x),\text{ for
			all }x\in \left[ \frac{j}{n-1},\frac{j+1}{n-1}\right] ,\text{ }j=0,1,2,...
	\end{equation*}%
	where $b_{n,k}\left( x\right) =\binom{n+k-1}{k}\frac{x^{k}}{\left(
		1+x\right) ^{n+k}}.$
\end{lem}

Let us define the following expression similar to \cite{max-prod Baskakov}.

For each $n\in \mathbb{N},n\geq 2$, $k,j\in \left\{ 0,1,2,...\right\} $ and $%
x\in \left[ \frac{j}{n-1},\frac{j+1}{n-1}\right] ,$ $x>0,$

\begin{eqnarray*}
	m_{k,n,j}\left( x\right) &:&=\dfrac{b_{n,k}\left( x\right) }{b_{n,j}\left(
		x\right) }=\frac{\binom{n+k-1}{k}}{\binom{n+j-1}{j}}\frac{x^{k}}{\left(
		1+x\right) ^{n+k}}\frac{\left( 1+x\right) ^{n+j}}{x^{j}} \\
	&=&\frac{\binom{n+k-1}{k}}{\binom{n+j-1}{j}}\left( \frac{x}{1+x}\right)
	^{k-j}.
\end{eqnarray*}%
And for $x=0$ let us denote $m_{0,n,0}\left( x\right) =1$ and $%
m_{k,n,0}\left( x\right) =0$ for all $k\in \left\{ 1,2,...\right\} .$

\begin{lem}
	\label{Lemma4.2} \cite{max-prod Baskakov} Let $n\in \mathbb{N},$ $n\geq 2.$\
	For all $k,$ $j\in \left\{ 0,1,2,...\right\} $ and $x\in \left[ \frac{j}{n-1}%
	,\frac{j+1}{n-1}\right] \ $we have%
	\begin{equation*}
		m_{k,n,j}(x)\leq 1.
	\end{equation*}
\end{lem}

\begin{rem}
	From Lemma \ref{Lemma3.1}, \ref{Lemma4.1}  and Lemma \ref{Lemma4.2} it is
	clear that $V_{n}^{(M)}(f)(x)$ satisfies for all $n\in \mathbb{N}$, $n\geq 2$
	all the hypothesis in Lemma \ref{Lemma2.1}, Corollary \ref{Corollary2.1} and
	Corollary \ref{Corollary2.2} for $I=\left[ 0,\infty \right) .$
\end{rem}

From Lemma \ref{Lemma3.1} we get, $V_{n}^{(M)}(f)(0)-f\left( 0\right) =0$
for all $n\geq 3,$ so in this part, we will consider $x>0$ in the notations,
proofs and statements of the all approximation results.

Again, let's define the following expressions similar to \cite{max-prod
	Baskakov}.

For each $n\in \mathbb{N},$ $n\geq 3$, $k,j\in \left\{ 0,1,2,...\right\} $
and $x\in \left[ \frac{j}{n-1},\frac{j+1}{n-1}\right] ,$

\begin{equation*}
	M_{k,n,j}\left( x\right) :=m_{k,n,j}\left( x\right) \left\vert \frac{k}{n}%
	-x\right\vert .
\end{equation*}%
It is clear that if $k\geq \frac{n}{n-1}\left( j+1\right) $ then we get 
\begin{equation*}
	M_{k,n,j}\left( x\right) =m_{k,n,j}\left( x\right) \left( \frac{k}{n}%
	-x\right)
\end{equation*}%
and if $k\leq \frac{n}{n-1}j$ then we have 
\begin{equation*}
	M_{k,n,j}\left( x\right) =m_{k,n,j}\left( x\right) \left( x-\frac{k}{n}%
	\right) .
\end{equation*}%
Also, for each $n\in \mathbb{N},$ $n\geq 3$, $k,j\in \mathbb{N}$, $k\geq 
\frac{n}{n-1}\left( j+1\right) $ and $x\in \left[ \frac{j}{n-1},\frac{j+1}{%
	n-1}\right] $%
\begin{equation*}
	\overline{M}_{k,n,j}\left( x\right) :=m_{k,n,j}\left( x\right) \left( \frac{k%
	}{n-1}-x\right)
\end{equation*}%
and for each $n\in \mathbb{N},$ $n\geq 3$, $k,j\in \mathbb{N}$, $k\leq \frac{%
	n}{n+1}j$ and $x\in \left[ \frac{j}{n-1},\frac{j+1}{n-1}\right] $%
\begin{equation*}
	\text{$\underline{M}$}_{k,n,j}\left( x\right) :=m_{k,n,j}\left( x\right)
	\left( x-\frac{k}{n-1}\right) .
\end{equation*}

Notice that, for all $k,j\in \left\{ 0,1,2,...\right\} $ and $n\geq 3$, the
above expressions were defined in \cite{max-prod Baskakov}.

At this point, let us recall the following lemma.

\begin{lem}
	\label{Lemma4.3} \cite{max-prod Baskakov} Let $x\in \left[ \frac{j}{n-1},%
	\frac{j+1}{n-1}\right] $ and $n\in \mathbb{N},$ $n\geq 3.$
	
	$(i)$ For all $k,j\in \left\{ 0,1,2,...\right\} $ with $k\geq \frac{n}{n-1}%
	\left( j+1\right) $ we have%
	\begin{equation*}
		M_{k,n,j}\left( x\right) \leq \overline{M}_{k,n,j}\left( x\right) .
	\end{equation*}%
	$(ii)$ For all $k,j\in \mathbb{N}$ with $k\geq \frac{n}{n-2}\left(
	j+1\right) $ we have%
	\begin{equation*}
		\overline{M}_{k,n,j}\left( x\right) \leq 2M_{k,n,j}\left( x\right) .
	\end{equation*}%
	$(iii)$ For all $k,j\in \mathbb{N}$ with $k\leq \frac{n}{n+1}j$ we have%
	\begin{equation*}
		\text{$\underline{M}$}_{k,n,j}\left( x\right) \leq M_{k,n,j}\left( x\right)
		\leq 2\text{$\underline{M}$}_{k,n,j}\left( x\right) .
	\end{equation*}
\end{lem}

Now, we will give our first main result of this part which is proven by not
only using the proof techniques given in \cite{max-prod Baskakov} but also
using the induction method. Especially, the second parts of $(i)$ and $(ii)$
in the proof are different from the proof techniques given in \cite{max-prod
	Baskakov}.

\begin{lem}
	\label{Lemma4.4} Let $x\in \left[ \frac{j}{n-1},\frac{j+1}{n-1}\right] $ and 
	$n\in \mathbb{N},$ $n\geq 3$ and$\ \alpha \in \left\{ 2,3,...\right\} .$
	
	$(i)$ If $j\in \left\{ 0,1,2,...\right\} $ is such that $k\geq \frac{n}{n-1}%
	\left( j+1\right) $ and 
	\begin{equation*}
		(k-j)^{\alpha }\geq \frac{\left( n+j\right) \left( k+1\right) }{\left(
			n-1\right) },
	\end{equation*}%
	then we have 
	\begin{equation*}
		\overline{M}_{k,n,j}(x)\geq \overline{M}_{k+1,n,j}(x).
	\end{equation*}%
	$(ii)$ If $k\in \left\{ 1,2,...,j\right\} $ is such that $k\leq \frac{n}{n+1}%
	j$ and 
	\begin{equation*}
		\left( j-k\right) ^{\alpha }\geq \frac{k\left( n+j-1\right) }{\left(
			n-1\right) },
	\end{equation*}%
	then we have 
	\begin{equation*}
		\text{$\underline{M}$}_{k,n,j}(x)\geq \text{$\underline{M}$}_{k-1,n,j}(x).
	\end{equation*}
\end{lem}

\begin{pf}
	$(i)$ From the case $(i)$ of Lemma 3.2 in \cite{max-prod Baskakov}, we can
	write 
	\begin{equation*}
		\dfrac{\overline{M}_{k,n,j}(x)}{\overline{M}_{k+1,n,j}(x)}\geq \frac{k+1}{n+k%
		}\dfrac{n+j}{j+1}\frac{k-j-1}{k-j}.
	\end{equation*}%
	After this point we will use a different proof technique from \cite{max-prod
		Baskakov}.
	
	By using the induction method, let's show that, the following inequality 
	\begin{equation}
		\frac{k+1}{n+k}\dfrac{n+j}{j+1}\frac{k-j-1}{k-j}\geq 1  \label{3.1}
	\end{equation}%
	holds for $(k-j)^{\alpha }\geq \frac{\left( n+j\right) \left( k+1\right) }{%
		\left( n-1\right) }.$
	
	For $\alpha =2$, since the condition $(k-j)^{2}\geq \frac{\left( n+j\right)
		\left( k+1\right) }{\left( n-1\right) }$ holds, because of%
	\begin{equation*}
		\begin{tabular}{l}
			$\left( k^{2}-2kj+j^{2}\right) \left( n-1\right) \geq \left( n+j\right)
			\left( k+1\right) \medskip $ \\ 
			$nk^{2}-2nkj+nj^{2}-k^{2}+kj-j^{2}-nk-n-j\geq 0\medskip $ \\ 
			
			$\left( nk+kj+n+j\right) \left( k-j-1\right) \geq \left( nj+n+kj+k\right)
			\left( k-j\right) .$%
		\end{tabular}%
	\end{equation*}%
	So we get $\left( k+1\right) \left( n+j\right) \left( k-j-1\right) \geq
	\left( n+k\right) \left( j+1\right) \left( k-j\right) $ therefore we obtain
	the inequality (\ref{3.1}) for $\alpha =2.$
	
	Now, we assume that the inequality (\ref{3.1}) is provided for $\alpha -1.$
	It follows $\frac{k+1}{n+k}\frac{n+j}{j+1}\frac{k-j-1}{k-j}\geq 1$ when $%
	(k-j)^{\alpha -1}\geq \frac{\left( n+j\right) \left( k+1\right) }{\left(
		n-1\right) }.$
	
	Since $k\geq \frac{n}{n-1}\left( j+1\right) $ we have $nk-k\geq nj+n$
	clearly it follows $n\left( k-j\right) \geq n+k$ or $\left( k-j\right) \geq
	1+\frac{k}{n}\geq 1,$ and $\alpha =2,3,...$ then $\left( k-j\right) ^{\alpha
	}\geq (k-j)^{\alpha -1}\geq \frac{\left( n+j\right) \left( k+1\right) }{%
		\left( n-1\right) }$ is true for $\alpha ,$ hence, for arbitrary $\alpha
	=2,3,...$ the inequality (\ref{3.1}) is provided when $(k-j)^{\alpha }\geq 
	\frac{\left( n+j\right) \left( k+1\right) }{\left( n-1\right) }.$
	
	So we obtain,%
	\begin{equation*}
		\dfrac{\overline{M}_{k,n,j}(x)}{\overline{M}_{k+1,n,j}(x)}\geq \frac{k+1}{n+k%
		}\dfrac{n+j}{j+1}\frac{k-j-1}{k-j}\geq 1.
	\end{equation*}%
	$(ii)$ From the case $(ii)$ of Lemma 3.2 in \cite{max-prod Baskakov}, we can
	write 
	\begin{equation*}
		\dfrac{\text{$\underline{M}$}_{k,n,j}(x)}{\text{$\underline{M}$}_{k-1,n,j}(x)%
		}\geq \frac{n+k-1}{k}\frac{j}{n+j-1}\frac{j-k}{j-k+1}.
	\end{equation*}%
	After this point we will use the our proof technique again.
	
	By using the induction method, let's show that the following inequality 
	\begin{equation}
		\frac{n+k-1}{k}\frac{j}{n+j-1}\frac{j-k}{j-k+1}\geq 1  \label{3.2}
	\end{equation}%
	holds for $\left( j-k\right) ^{\alpha }\geq \frac{k\left( n+j-1\right) }{%
		\left( n-1\right) }=k\left( \frac{j}{n-1}+1\right) $.
	
	For $\alpha =2$, because of the condition $\left( j-k\right) ^{2}\geq \frac{%
		k\left( n+j-1\right) }{\left( n-1\right) }$, i.e $\left( j-k\right)
	^{2}\left( n-1\right) -k\left( n+j-1\right) \geq 0.$ After simple calculus
	we get%
	\begin{equation*}
		\begin{tabular}{l}
			$nj^{2}-2nkj+nk^{2}+2kj-j^{2}-k^{2}-nk-kj+k\geq 0\medskip $ \\ 
		
			$\left( n+k-1\right) \left( j^{2}-kj\right) \geq \left( nk+kj-k\right)
			\left( j-k+1\right) \medskip $ \\ 
			$\left( n+k-1\right) j\left( j-k\right) \geq k\left( n+j-1\right) \left(
			j-k+1\right) $%
		\end{tabular}%
	\end{equation*}%
	so we see that (\ref{3.2}) is satisfied.
	
	Now, we assume that (\ref{3.2}) is correct for $\alpha -1.$ Hence, $\frac{%
		n+k-1}{k}\frac{j}{n+j-1}\frac{j-k}{j-k+1}\geq 1$ is provided when $\left(
	j-k\right) ^{\alpha -1}\geq \frac{k\left( n+j-1\right) }{\left( n-1\right) }$
	and $k\leq \frac{n}{n+1}j.$ We can assume that $k\leq j-1$, because when $%
	k=j,$ other conditions, which are $k\leq \frac{n}{n+1}j$ and $\left(
	j-k\right) ^{\alpha }\geq \frac{k\left( n+j-1\right) }{\left( n-1\right) },$
	are not satisfied.
	
	So, since $k\leq j-1$ or $1\leq j-k$ it follows that $\left( j-k\right)
	^{\alpha -1}\leq \left( j-k\right) ^{\alpha }\geq \frac{k\left( n+j-1\right) 
	}{\left( n-1\right) }$ is true for $\alpha =2,3,...$ then the desired
	inequality is provided for $\left( j-k\right) ^{\alpha }\geq \frac{k\left(
		n+j-1\right) }{\left( n-1\right) }.$ So we obtain,%
	\begin{equation*}
		\dfrac{\text{$\underline{M}$}_{k,n,j}(x)}{\text{$\underline{M}$}_{k-1,n,j}(x)%
		}\geq \frac{n+k-1}{k}\frac{j}{n+j-1}\frac{j-k}{j-k+1}\geq 1,
	\end{equation*}%
	which gives the desired result.
\end{pf}

\section{Pointwise Rate of Convergence}

Let us take a $x_{0}$ fixed point on the interval $[0,\infty )$. The main
aim of this section is to obtain a better order of pointwise approximation
for the operators $V_{n}^{(M)}(f)(x_{0})$ to the function $f(x_{0})$ by
means of the classical modulus of continuity. According to the following
theorem we can say that the order of pointwise approximation can be improved
when the $\alpha $ is big enough. Moreover if we choose as $\alpha =2$,
these approximation results turn out to be the results in \cite{max-prod
	Baskakov}.

\begin{thm}
	\label{Theorem5.1} Let $f:\left[ 0,\infty \right) \rightarrow \mathbb{R}_{+}$
	be bounded and continuous. Then for any fixed point $x_{0}$ on the interval $%
	[0,\infty ),$ which also satisfy $x_{0}^{\alpha -2}\leq n-1$ and $n\geq
	j^{\alpha -1},$ $n\geq 4$ we have the following order of approximation for
	the operators (\ref{2.1}) to the function $f$ by means of the modulus of
	continuity:%
	\begin{equation*}
		\left\vert V_{n}^{(M)}(f)(x_{0})-f(x_{0})\right\vert \leq (1+6\left[
		x_{0}\left( 1+x_{0}\right) \right] ^{\frac{1}{\alpha }})\text{ }\omega
		\left( f;\frac{1}{\left( n-1\right) ^{1-\frac{1}{\alpha }}}\right) ,
	\end{equation*}%
		for all $n\in \mathbb{N},\text{ }n\geq 4, $ where $\omega \left( f;\delta \right) $ is the classical modulus of
	continuity defined by (\ref{1.1}) and $\alpha =2,3,...$ .
\end{thm}

\begin{pf}
	Since nonlinear max-product Baskakov operators satisfy the conditions in
	Corollary \ref{Corollary2.2}, for any $x_{0}\in $ $[0,\infty ),$ using the
	properties of $\omega \left( f;\delta \right) ,$ we get%
	\begin{equation}
		\left\vert V_{n}^{(M)}\left( f\right) \left( x_{0}\right) -f\left(
		x_{0}\right) \right\vert \leq \left[ 1+\frac{1}{\delta _{n}}%
		V_{n}^{(M)}\left( \varphi _{x_{0}}\right) \left( x_{0}\right) \right]
		\,\omega \left( f,\delta \right) ,  \label{4.1}
	\end{equation}%
	where $\varphi _{x_{0}}\left( t\right) =\left\vert t-x_{0}\right\vert .$ At
	this point let us denote%
	\begin{equation*}
		E_{n}\left( x_{0}\right) :=V_{n}^{(M)}\left( \varphi _{x_{0}}\right) \left(
		x_{0}\right) =\frac{\bigvee\limits_{k=0}^{\infty }b_{n,k}\left( x_{0}\right)
			\left\vert \frac{k}{n}-x_{0}\right\vert }{\bigvee\limits_{k=0}^{\infty
			}b_{n,k}\left( x_{0}\right) },\text{ }x_{0}\in \left[ 0,\infty \right) .
	\end{equation*}%
	Let $x_{0}\in \left[ \frac{j}{n-1},\frac{j+1}{n-1}\right] ,$ where $j\in
	\left\{ 0,1,...\right\} $ is fixed, arbitary. By Lemma \ref{Lemma4.1} we
	easily obtain 
	\begin{equation*}
		E_{n}\left( x_{0}\right) =\max\limits_{k=0,1,...}\left\{ M_{k,n,j}\left(
		x_{0}\right) \right\} .
	\end{equation*}%
	Firstly let's check for $j=0$, $M_{k,n,0}\left( x_{0}\right) =\frac{\binom{%
			n+k-1}{k}}{\binom{n+0-1}{0}}\left( \frac{x_{0}}{1+x_{0}}\right)
	^{k-0}\left\vert \frac{k}{n}-x_{0}\right\vert $\ where $x_{0}\in \left[ 0,%
	\frac{1}{n-1}\right] $ and $\alpha =2,3,...$ .
	
	For $k=0$, we get%
	\begin{equation*}
		M_{0,n,0}\left( x_{0}\right) =x_{0}=x_{0}^{\frac{1}{\alpha }}x_{0}^{1-\frac{1%
			}{\alpha }}\leq \frac{x_{0}^{\frac{1}{\alpha }}}{\left( n-1\right) ^{1-\frac{%
					1}{\alpha }}}.
	\end{equation*}%
	And for $k=1,$ we get 
	\begin{equation*}
		M_{1,n,0}\left( x_{0}\right) =\binom{n}{1}\left( \frac{x_{0}}{1+x_{0}}%
		\right) \left\vert \frac{1}{n}-x_{0}\right\vert \leq n\frac{x_{0}}{1+x_{0}}%
		\frac{1}{n}\leq x_{0}\leq \frac{x_{0}^{\frac{1}{\alpha }}}{\left( n-1\right)
			^{1-\frac{1}{\alpha }}}.
	\end{equation*}%
	Now suppose that $k\geq 2.$ For $j=0$ and $k\geq 2$, we see that all of the
	Lemma \ref{Lemma4.3} $(i)$'s hypotheses are satisfied. Thus we get $%
	M_{k,n,0}\left( x_{0}\right) \leq \overline{M}_{k,n,0}(x).$ Also by Lemma %
	\ref{Lemma4.4} $(i)$, for $j=0$ it follows that $\overline{M}_{k,n,0}(x)\geq 
	\overline{M}_{k+1,n,0}(x)$ for every $k\geq 2$ such that $k^{\alpha }\geq 
	\frac{n\left( k+1\right) }{n-1}$ then we have $\left( n-1\right) k^{\alpha
	}-nk-n\geq 0,$ when $\alpha =2,3,...$ . Since the function $f_{n}\left(
	t\right) =\left( n-1\right) t^{\alpha }-nt-n,$ $t\geq 1$ is nondecreasing,
	really since $f_{n}^{\prime }\left( t\right) =\alpha \left( n-1\right)
	t^{\alpha -1}-n\geq 0,$ and because%
	\begin{equation*}
		f_{n}\left( n^{\frac{1}{\alpha }}\right) =\left( n-1\right) n-n\,n^{\frac{1}{%
				\alpha }}-n=n\left( n-2-n^{\frac{1}{\alpha }}\right) \geq 0,\text{ }n\geq 4,
	\end{equation*}%
	thus, it follows $\overline{M}_{k,n,0}(x)\geq \overline{M}_{k+1,n,0}(x)$ for
	every $k\in \mathbb{N}$, $k\geq n^{\frac{1}{\alpha }}.$ Let us denote 
	\begin{equation*}
		A=\left\{ k\in \mathbb{N},\text{ }2\leq k\leq n^{\frac{1}{\alpha }}+1\right\}
	\end{equation*}%
	and let $k\in A.$ Since $4\leq n$ then $0\leq k\left( n-3\right) $ it
	follows that $2nk\leq 3nk-3k$ or $\frac{k}{n-1}\leq \frac{3k}{2n},$ and by
	Lemma \ref{Lemma4.2} then we obtain%
	\begin{eqnarray*}
		\overline{M}_{k,n,0}(x_{0}) &=&\binom{n+k-1}{k}\left( \frac{x_{0}}{1+x_{0}}%
		\right) ^{k}\left( \frac{k}{n-1}-x_{0}\right) \\
		&\leq &\binom{n+k-1}{k}\left( \frac{x_{0}}{1+x_{0}}\right) ^{k}\frac{k}{n-1}
		\\
		&\leq &\binom{n+k-1}{k}\left( \frac{x_{0}}{1+x_{0}}\right) ^{k}\frac{3k}{2n}
		\\
		&=&\frac{\left( n+k-1\right) !}{k!\left( n-1\right) !}\frac{3k}{2n}\left( 
		\frac{x_{0}}{1+x_{0}}\right) ^{k} \\
		&=&\frac{\left( n+k-1\right) !}{\left( k-1\right) !n!}\frac{3}{2}\left( 
		\frac{x_{0}}{1+x_{0}}\right) ^{k} \\
		&=&\frac{3}{2}\binom{n+k-1}{k-1}\left( \frac{x_{0}}{1+x_{0}}\right) ^{k-1}%
		\frac{x_{0}}{1+x_{0}} \\
		&=&\frac{3}{2}\binom{n+k-1}{k-1}\left( \frac{\frac{1}{n}}{1+\frac{1}{n}}%
		\right) ^{k-1}\left( \frac{x_{0}}{1+x_{0}}\frac{1+\frac{1}{n}}{\frac{1}{n}}%
		\right) ^{k-1}\frac{x_{0}}{1+x_{0}} \\
		&=&\frac{3}{2}m_{k-1,n+1,0}\left( \frac{1}{n}\right) \left( \frac{\left(
			n+1\right) x_{0}}{1+x_{0}}\right) ^{k-1}\frac{x_{0}}{1+x_{0}} \\
		&\leq &\frac{3}{2}\left( \frac{\left( n+1\right) x_{0}}{1+x_{0}}\right)
		^{k-1}\frac{x_{0}}{\left( 1+x_{0}\right) },
	\end{eqnarray*}%
	taking into account that $\frac{1}{n+1}=\frac{\frac{1}{n}}{1+\frac{1}{n}}$.
	If we deneote%
	\begin{equation*}
		g_{n,k}\left( x_{0}\right) =\left( \frac{\left( n+1\right) x_{0}}{1+x_{0}}%
		\right) ^{k-1},
	\end{equation*}%
	then we see that the function $g_{n,k}\left( x_{0}\right) $ is nondecreasing
	on the interval $\left[ 0,\frac{1}{n-1}\right] $. Really since 
	\begin{eqnarray*}
		g_{n,k}^{\prime }\left( x_{0}\right) &=&\left( k-1\right) \left( \frac{%
			\left( n+1\right) x_{0}}{1+x_{0}}\right) ^{k-2}\left( \frac{\left(
			n+1\right) \left( 1+x_{0}\right) -\left( n+1\right) x_{0}}{\left(
			1+x_{0}\right) ^{2}}\right) \\
		&=&\left( k-1\right) \left( \frac{\left( n+1\right) x_{0}}{1+x_{0}}\right)
		^{k-2}\frac{n+1}{\left( 1+x_{0}\right) ^{2}}\geq 0.
	\end{eqnarray*}%
	Using this property, and since $x_{0}\leq \frac{1}{n-1}$ we obtain%
	\begin{equation*}
		g_{n,k}\left( x_{0}\right) \leq g_{n,k}\left( \frac{1}{n-1}\right) =\left( 
		\frac{\frac{n+1}{n-1}}{1+\frac{1}{n-1}}\right) ^{k-1}=\left( \frac{n+1}{n}%
		\right) ^{k-1},
	\end{equation*}%
	for all $x_{0}\in \left[ 0,\frac{1}{n-1}\right] .$ Then 
	\begin{eqnarray*}
		\overline{M}_{k,n,0}(x_{0}) &\leq &\frac{3}{2}\left( \frac{n+1}{n}\right)
		^{k-1}\frac{x_{0}}{\left( 1+x_{0}\right) } \\
		&<&\frac{3}{2}\left( \frac{n+1}{n}\right) ^{n^{\frac{1}{\alpha }}}\frac{x_{0}%
		}{\left( 1+x_{0}\right) } \\
		&\leq &\frac{3}{2}\left( \frac{n+1}{n}\right) ^{n}\frac{x_{0}}{\left(
			1+x_{0}\right) } \\
		&<&\frac{3}{2}e\frac{x_{0}}{1+x_{0}}<\frac{3e}{2}x_{0} \\
		&\leq &5\frac{x_{0}^{\frac{1}{\alpha }}}{\left( n-1\right) ^{1-\frac{1}{%
					\alpha }}}
	\end{eqnarray*}%
	taking into account that $k+1\leq n^{\frac{1}{\alpha }}$ and $%
	\lim_{n\rightarrow \infty }\left( 1+\frac{1}{n}\right) ^{n}=e$ and $\frac{3}{%
		2}e<5.$
	
	So, we find an upper estimate for any $k=0,1,2,...$ 
	\begin{eqnarray*}
		E_{n}\left( x_{0}\right) &=&\max\limits_{k=0,1,...}\left\{ M_{k,n,0}\left(
		x_{0}\right) \right\} \\
		&\leq &\max \left\{ M_{0,n,0}\left( x_{0}\right) ,M_{1,n,0}\left(
		x_{0}\right) ,\max\limits_{k=2,3,...}\left\{ \overline{M}_{k,n,0}\left(
		x_{0}\right) \right\} \right\} \\
		&=&\max \left\{ M_{0,n,0}\left( x_{0}\right) ,M_{1,n,0}\left( x_{0}\right)
		,\max\limits_{k\in A}\left\{ \overline{M}_{k,n,0}\left( x_{0}\right)
		\right\} \right\} \\
		&<&5\frac{x_{0}^{\frac{1}{\alpha }}}{\left( n-1\right) ^{1-\frac{1}{\alpha }}%
		}\leq 5\frac{\left[ x_{0}\left( 1+x_{0}\right) \right] ^{\frac{1}{\alpha }}}{%
			\left( n-1\right) ^{1-\frac{1}{\alpha }}},
	\end{eqnarray*}%
	when $j=0.$
	
	As a result, it remains to find an upper estimate for each $M_{k,n,j}\left(
	x_{0}\right) $ when $j=1,2,...$ is fixed, $x_{0}\in \left[ \frac{j}{n-1},%
	\frac{j+1}{n-1}\right] $, $k\in \left\{ 0,1,...\right\} $ and $\alpha
	=2,3,...$ .
	
	In fact we will show that 
	\begin{equation}
		M_{k,n,j}\left( x_{0}\right) \leq 6\frac{\left[ x_{0}\left( 1+x_{0}\right) %
			\right] ^{\frac{1}{\alpha }}}{\left( n-1\right) ^{1-\frac{1}{\alpha }}}
		\label{4.2}
	\end{equation}%
	for all $x\in \left[ \frac{j}{n-1},\frac{j+1}{n-1}\right] ,$ $k=0,1,2,...$
	it implies that directly 
	\begin{equation*}
		E_{n}\left( x_{0}\right) \leq 6\frac{\left[ x_{0}\left( 1+x_{0}\right) %
			\right] ^{\frac{1}{\alpha }}}{\left( n-1\right) ^{1-\frac{1}{\alpha }}},%
		\text{ for all }x\in \left[ 0,\infty \right) \text{ and }n\in \mathbb{N}%
		\text{, }n\geq 4,
	\end{equation*}%
	and taking $\delta _{n}=\frac{1}{\left( n-1\right) ^{1-\frac{1}{\alpha }}}$
	in (\ref{4.1}),\ we obtain the estimate in the statement immediately.
	
	So, in order to completing the prove of (\ref{4.2}), we consider the
	following cases:
	
	$1)$ $\frac{n}{n+1}j\leq k\leq \frac{n}{n-1}\left( j+1\right) ;$
	
	$2)$ $k>\frac{n}{n-1}\left( j+1\right) ;$
	
	$3)$ $k<\frac{n}{n+1}j.$
	
	Case $1).$ We get 
	\begin{equation*}
		\frac{k}{n}-x_{0}\leq \frac{\frac{n}{n-1}\left( j+1\right) }{n}-x_{0}\leq 
		\frac{j+1}{n-1}-\frac{j}{n-1}=\frac{1}{n-1}\leq \frac{2x_{0}}{n-1}+\frac{1}{%
			n-1}.
	\end{equation*}%
	On the other hand 
	\begin{eqnarray*}
		\frac{k}{n}-x_{0} &\geq &\frac{\frac{n}{n+1}j}{n}-x_{0}\geq \frac{j}{n+1}-%
		\frac{j+1}{n-1} \\
		&=&\frac{nj-j-nj-n-j-1}{\left( n-1\right) \left( n+1\right) }=\frac{-2j-n-1}{%
			\left( n-1\right) \left( n+1\right) } \\
		&=&\frac{-2j}{\left( n-1\right) \left( n+1\right) }-\frac{1}{n-1}\geq -\frac{%
			2x_{0}}{n+1}-\frac{1}{n-1} \\
		&\geq &-\frac{2x_{0}}{n-1}-\frac{1}{n-1}.
	\end{eqnarray*}%
	As a result we obtain 
	\begin{equation*}
		\left\vert \frac{k}{n}-x_{0}\right\vert \leq \frac{2x_{0}}{n-1}+\frac{1}{n-1}%
		.
	\end{equation*}%
	Since $x_{0}^{\alpha -2}\leq n-1$ from the hypothesis, and $n\geq 4$ we
	obtain $\dfrac{x_{0}}{n-1}\leq $ $\dfrac{\left[ x_{0}\left( 1+x_{0}\right) %
		\right] ^{\frac{1}{\alpha }}}{\left( n-1\right) ^{1-\frac{1}{\alpha }}}$ for
	all $x>0$ and $\alpha =2,3,...$.
	
	And also we get, 
	\begin{eqnarray*}
		\frac{1}{n-1} &=&\left( \frac{1}{n-1}\right) ^{\frac{1}{\alpha }}\left( 
		\frac{1}{n-1}\right) ^{1-\frac{1}{\alpha }}\leq \left( \frac{j}{n-1}\right)
		^{\frac{1}{\alpha }}\frac{1}{\left( n-1\right) ^{1-\frac{1}{\alpha }}} \\
		&\leq &\frac{\left( x_{0}\right) ^{\frac{1}{\alpha }}}{\left( n-1\right) ^{1-%
				\frac{1}{\alpha }}}\leq \dfrac{\left[ x_{0}\left( 1+x_{0}\right) \right] ^{%
				\frac{1}{\alpha }}}{\left( n-1\right) ^{1-\frac{1}{\alpha }}}.
	\end{eqnarray*}%
	Hence, it follows 
	\begin{equation*}
		M_{k,n,j}\left( x\right) =m_{k,n,j}\left( x\right) \left\vert \frac{k}{n}%
		-x_{0}\right\vert \leq 3\frac{\left[ x_{0}\left( 1+x_{0}\right) \right] ^{%
				\frac{1}{\alpha }}}{\left( n-1\right) ^{1-\frac{1}{\alpha }}}.
	\end{equation*}%
	Case $2).$ Subcase $a).$ Assume first that $\left( k-j\right) ^{\alpha }<%
	\frac{\left( n+j\right) \left( k+1\right) }{\left( n-1\right) }.$ If we
	denoting $k=j+\beta $, where $\beta \geq 1,$ the condition becomes%
	\begin{equation*}
		\begin{tabular}{l}
			$\beta ^{\alpha }<\frac{\left( n+j\right) \left( j+\beta +1\right) }{\left(
				n-1\right) }\medskip $ \\ 
			$\beta ^{\alpha }\left( n-1\right) -\left( n+j\right) \left( j+\beta
			+1\right) <0\medskip $ \\ 
			$\beta ^{\alpha }\left( n-1\right) -\beta \left( n+j\right) -\left(
			n+j\right) \left( j+1\right) <0.$%
		\end{tabular}%
	\end{equation*}%
	Let us define the function $f\left( t\right) =t^{\alpha }\left( n-1\right)
	-t\left( n+j\right) -\left( n+j\right) \left( j+1\right) ,$ $t\in \mathbb{R}$%
	. We claim that $f\left( \left[ \frac{3\left( j+1\right) \left( n+j\right) }{%
		n-1}\right] ^{\frac{1}{\alpha }}\right) >0$, which will imply $k-j=\beta <%
	\left[ \frac{3\left( j+1\right) \left( n+j\right) }{n-1}\right] ^{\frac{1}{%
			\alpha }}.$ After simple calculation we have%
	\begin{eqnarray*}
		&&f\left( \left[ \frac{3\left( j+1\right) \left( n+j\right) }{n-1}\right] ^{%
			\frac{1}{\alpha }}\right) \\
		&=&\frac{3\left( j+1\right) \left( n+j\right) }{n-1}(n-1)-\left[ \frac{%
			3\left( j+1\right) \left( n+j\right) }{n-1}\right] ^{\frac{1}{\alpha }%
		}\left( n+j\right) -\left( n+j\right) \left( j+1\right) \\
		&=&2\left( n+j\right) \left( j+1\right) -\left[ \frac{3\left( j+1\right)
			\left( n+j\right) }{n-1}\right] ^{\frac{1}{\alpha }}\left( n+j\right) \\
		&=&\left( n+j\right) \left( j+1\right) ^{\frac{1}{\alpha }}\left\{ 2\left(
		j+1\right) ^{1-\frac{1}{\alpha }}-\left( \frac{3\left( n+j\right) }{n-1}%
		\right) ^{\frac{1}{\alpha }}\right\} \\
		&=&\left( n+j\right) \left( j+1\right) ^{\frac{1}{\alpha }}\left\{ 2\left(
		j+1\right) ^{1-\frac{1}{\alpha }}-\left( 3+\frac{3j+3}{n-1}\right) ^{\frac{1%
			}{\alpha }}\right\} \\
		&\geq &\left( n+j\right) \left( j+1\right) ^{\frac{1}{\alpha }}\left\{
		2\left( j+1\right) ^{1-\frac{1}{\alpha }}-\left( 3+\frac{3j+3}{2}\right) ^{%
			\frac{1}{\alpha }}\right\} \\
		&=&\left( n+j\right) \left( j+1\right) ^{\frac{1}{\alpha }}\left\{ 2\left(
		j+1\right) ^{1-\frac{1}{\alpha }}-\left( \frac{3j+9}{2}\right) ^{\frac{1}{%
				\alpha }}\right\} >0,
	\end{eqnarray*}%
	taking into account that this conclusion. For all $\alpha \geq 2$ and $j\geq
	1$ it is clear $2^{\alpha +1}\left( j+1\right) ^{\alpha -1}>3j+9.$ It
	follows that $2^{\alpha }\left( j+1\right) ^{\alpha -1}>\frac{3j+9}{2}$ and
	we obtain $2\left( j+1\right) ^{1-\frac{1}{\alpha }}>\left( \frac{3j+9}{2}%
	\right) ^{\frac{1}{\alpha }}.$
	
	Based on the findings above, we get 
	\begin{eqnarray*}
		\overline{M}_{k,n,j}\left( x_{0}\right) &=&m_{k,n,j}\left( x_{0}\right)
		\left( \frac{k}{n-1}-x_{0}\right) \leq \frac{k}{n-1}-x_{0} \\
		&\leq &\frac{k}{n-1}-\frac{j}{n-1}=\frac{k-j}{n-1}=\frac{\beta }{n-1} \\
		&<&\frac{\left[ \frac{3\left( j+1\right) \left( n+j\right) }{n-1}\right] ^{%
				\frac{1}{\alpha }}}{n-1}=\frac{\left[ 3\left( j+1\right) \left( n+j\right) %
			\right] ^{\frac{1}{\alpha }}}{\left( n-1\right) ^{1+\frac{1}{\alpha }}}\leq 
		\frac{\left[ 6j\left( n+j\right) \right] ^{\frac{1}{\alpha }}}{\left(
			n-1\right) ^{1+\frac{1}{\alpha }}} \\
		&=&\frac{1}{\left( n-1\right) ^{1-\frac{1}{\alpha }}}\left( \frac{6j}{n-1}%
		\right) ^{\frac{1}{\alpha }}\left( \frac{n+j}{n-1}\right) ^{\frac{1}{\alpha }%
		} \\
		&=&\frac{1}{\left( n-1\right) ^{1-\frac{1}{\alpha }}}\left( \frac{6j}{n-1}%
		\right) ^{\frac{1}{\alpha }}\left( \frac{n+j-1}{n-1}\right) ^{\frac{1}{%
				\alpha }}\left( \frac{n+j}{n+j-1}\right) ^{\frac{1}{\alpha }} \\
		&\leq &\frac{1}{\left( n-1\right) ^{1-\frac{1}{\alpha }}}\left(
		6x_{0}\right) ^{\frac{1}{\alpha }}\left( 1+x_{0}\right) ^{\frac{1}{\alpha }%
		}\left( \frac{4}{3}\right) ^{\frac{1}{\alpha }} \\
		&=&2^{\frac{3}{\alpha }}\frac{\left[ x_{0}\left( 1+x_{0}\right) \right] ^{%
				\frac{1}{\alpha }}}{\left( n-1\right) ^{1-\frac{1}{\alpha }}},
	\end{eqnarray*}%
	we used that $\frac{n+j-1}{n-1}=1+\frac{j}{n-1}\leq 1+x_{0}$ and $n\geq 3,$ $%
	j\geq 1$ so $4\leq n+j$ in the above inequality.
	
	Subcase $b).$ Assume now that $\left( k-j\right) ^{\alpha }\geq \frac{\left(
		n+j\right) \left( k+1\right) }{\left( n-1\right) },$ it means $\left(
	k-j\right) ^{\alpha }\left( n-1\right) -\left( n+j\right) \left( k+1\right)
	\geq 0.$ Because $n$ and $j$ are fixed, we can define the real function 
	\begin{equation*}
		g\left( t\right) :=\left( t-j\right) ^{\alpha }\left( n-1\right) -\left(
		n+j\right) \left( t+1\right) ,
	\end{equation*}%
	for all $t\in \mathbb{R}$. For $t\geq \frac{n}{n-1}\left( j+1\right) $, $%
	g\left( t\right) $ is nondecreasing on the interval $\left[ \frac{n}{n-1}%
	\left( j+1\right) ,\infty \right) .$ Really since, 
	\begin{eqnarray*}
		g^{\prime }\left( t\right) &=&\alpha \left( t-j\right) ^{\alpha -1}\left(
		n-1\right) -n-j \\
		&\geq &\alpha \left( \frac{n}{n-1}\left( j+1\right) -j\right) ^{\alpha
			-1}\left( n-1\right) -n-j \\
		&=&\alpha \left( \frac{n+j}{n-1}\right) ^{\alpha -1}\left( n-1\right) -n-j \\
		&=&\left( n+j\right) \left[ \alpha \left( \frac{n+j}{n-1}\right) ^{\alpha
			-2}-1\right] >0.
	\end{eqnarray*}%
	In the inequality above, $\alpha \geq 2$ i.e $\frac{1}{\alpha }<1,$ and
	since $n+j>n-1$ or $\frac{n+j}{n-1}>1$ then we have $\left( \frac{n+j}{n-1}%
	\right) ^{\alpha -2}>\frac{1}{\alpha },$\ is taken into consideration.
	
	And since $\lim\limits_{t\rightarrow \infty }g\left( t\right) =\infty ,$ by
	the monotonicity of $g$ too, it follows there exists $\bar{k}\in \mathbb{N}$%
	, $\bar{k}>\frac{n}{n-1}\left( j+1\right) $ of minimum value, satisfying the
	inequality $g\left( \bar{k}\right) =\left( \bar{k}-j\right) ^{\alpha }\left(
	n-1\right) -\left( n+j\right) \left( \bar{k}+1\right) \geq 0.$ Denote $k_{1}=%
	\bar{k}+1,$ where evidently $k_{1}\geq j+1.$ If $k_{1}\geq \frac{n}{n-1}%
	\left( j+1\right) ,$ then from the properties of $g$ and by the way we
	choose $\bar{k}$ it results that $g\left( k_{1}\right) <0.$ If $k_{1}<\frac{n%
	}{n-1}\left( j+1\right) ,$ then $j<k_{1}<\frac{n}{n-1}\left( j+1\right) .$
	Since $g$ is a polynomial function and because $g\left( j\right) <0$ and $%
	g\left( \frac{n}{n-1}\left( j+1\right) \right) <0$, we immediately obtain
	the same conclusion as in the previous case, which is $g\left( k_{1}\right)
	<0$ or equivalently $\beta ^{\alpha }\left( n-1\right) -\left( n+j\right)
	\left( j+\beta +1\right) <0$, where $k_{1}=j+\beta .$ Using the same method
	as in subcase $a)$ we have $k_{1}-j<\left[ \frac{3\left( j+1\right) \left(
		n+j\right) }{n-1}\right] ^{\frac{1}{\alpha }}.$ Then%
	\begin{eqnarray*}
		\overline{M}_{\bar{k},n,j}\left( x_{0}\right) &=&m_{\bar{k},n,j}\left(
		x_{0}\right) \left( \frac{\bar{k}}{n-1}-x_{0}\right) \leq \frac{\bar{k}}{n-1}%
		-x_{0} \\
		&\leq &\frac{\bar{k}}{n-1}-\frac{j}{n-1}=\frac{\bar{k}-j}{n-1}=\frac{k_{1}-j%
		}{n-1}+\frac{1}{n-1} \\
		&\leq &2^{\frac{3}{\alpha }}\frac{\left[ x_{0}\left( 1+x_{0}\right) \right]
			^{\frac{1}{\alpha }}}{\left( n-1\right) ^{1-\frac{1}{\alpha }}}+\frac{1}{n-1}
		\\
		&\leq &2^{\frac{3}{\alpha }}\frac{\left[ x_{0}\left( 1+x_{0}\right) \right]
			^{\frac{1}{\alpha }}}{\left( n-1\right) ^{1-\frac{1}{\alpha }}}+\frac{\left[
			x_{0}\left( 1+x_{0}\right) \right] ^{\frac{1}{\alpha }}}{\left( n-1\right)
			^{1-\frac{1}{\alpha }}} \\
		&\leq &4\frac{\left[ x_{0}\left( 1+x_{0}\right) \right] ^{\frac{1}{\alpha }}%
		}{\left( n-1\right) ^{1-\frac{1}{\alpha }}},
	\end{eqnarray*}%
	taking into account that, since $\alpha \geq 2,$ then $2^{\frac{3}{\alpha }%
	}+1\leq 4.$ By the Lemma \ref{Lemma4.4} $(i)$ it follows that $\overline{M}_{%
		\bar{k},n,j}\left( x_{0}\right) \geq \overline{M}_{\bar{k}+1,n,j}\left(
	x_{0}\right) \geq ...$ . Therefore we obtain $\overline{M}_{\bar{k}%
		,n,j}\left( x_{0}\right) <4\frac{\left[ x_{0}\left( 1+x_{0}\right) \right] ^{%
			\frac{1}{\alpha }}}{\left( n-1\right) ^{1-\frac{1}{\alpha }}}$ for any $k\in
	\left\{ \bar{k},\bar{k}+1,...\right\} .$
	
	As a result, in both subcases, by Lemma \ref{Lemma4.3} $(i)$ too, we have%
	\begin{equation*}
		M_{k,n,j}\left( x_{0}\right) <4\frac{\left[ x_{0}\left( 1+x_{0}\right) %
			\right] ^{\frac{1}{\alpha }}}{\left( n-1\right) ^{1-\frac{1}{\alpha }}}.
	\end{equation*}%
	Case $3).$ Subcase $a).$ Assume first that $\left( j-k\right) ^{\alpha }<%
	\frac{k\left( n+j-1\right) }{\left( n-1\right) }.$ If we denoting $k=j-\beta 
	$, where $\beta \geq 1,$ the condition becomes 
	\begin{equation*}
		\begin{tabular}{l}
			$\beta ^{\alpha }<\frac{\left( j-\beta \right) \left( n+j-1\right) }{\left(
				n-1\right) }\medskip $ \\ 
			$\beta ^{\alpha }\left( n-1\right) -\left( j-\beta \right) \left(
			n+j-1\right) <0.$%
		\end{tabular}%
	\end{equation*}%
	Let us define the function $f\left( t\right) =t^{\alpha }\left( n-1\right)
	-\left( j-t\right) \left( n+j-1\right) ,$ $t\in \mathbb{R}$. We claim that $%
	f\left( \left[ \frac{j\left( n+j-1\right) }{n-1}\right] ^{\frac{1}{\alpha }%
	}\right) >0$ which will imply $j-k=\beta <\left[ \frac{j\left( n+j-1\right) 
	}{n-1}\right] ^{\frac{1}{\alpha }}.$ After simple calculation we have%
	\begin{eqnarray*}
		&&f\left( \left[ \frac{j\left( n+j-1\right) }{n-1}\right] ^{\frac{1}{\alpha }%
		}\right) \\
		&=&\frac{j\left( n+j-1\right) }{n-1}(n-1)-\left( j-\left[ \frac{j\left(
			n+j-1\right) }{n-1}\right] ^{\frac{1}{\alpha }}\right) \left( n+j-1\right) \\
		&=&j\left( n+j-1\right) -j\left( n+j-1\right) +\left[ \frac{j\left(
			n+j-1\right) }{n-1}\right] ^{\frac{1}{\alpha }}\left( n+j-1\right) \\
		&=&\left[ \frac{j\left( n+j-1\right) }{n-1}\right] ^{\frac{1}{\alpha }%
		}\left( n+j-1\right) >0.
	\end{eqnarray*}%
	Then we obtain 
	\begin{eqnarray*}
		\text{$\underline{M}$}_{k,n,j}\left( x_{0}\right) &=&m_{k,n,j}\left(
		x_{0}\right) \left( x_{0}-\frac{k}{n-1}\right) \leq \frac{j+1}{n-1}-\frac{k}{%
			n-1} \\
		&=&\frac{j-k}{n-1}+\frac{1}{n-1}<\frac{\left[ \frac{j\left( n+j-1\right) }{%
				n-1}\right] ^{\frac{1}{\alpha }}}{n-1}+\frac{1}{n-1} \\
		&=&\frac{1}{\left( n-1\right) ^{1-\frac{1}{\alpha }}}\left( \frac{j\left(
			n+j-1\right) }{n-1}\right) ^{\frac{1}{\alpha }}\left( \frac{1}{n-1}\right) ^{%
			\frac{1}{\alpha }}+\frac{1}{n-1} \\
		&=&\frac{1}{\left( n-1\right) ^{1-\frac{1}{\alpha }}}\left( \frac{j}{n-1}%
		\right) ^{\frac{1}{\alpha }}\left( 1+\frac{j}{n-1}\right) ^{\frac{1}{\alpha }%
		}+\frac{1}{n-1} \\
		&\leq &\frac{1}{\left( n-1\right) ^{1-\frac{1}{\alpha }}}\left( x_{0}\right)
		^{\frac{1}{\alpha }}\left( 1+x_{0}\right) ^{\frac{1}{\alpha }}+\frac{\left[
			x_{0}\left( 1+x_{0}\right) \right] ^{\frac{1}{\alpha }}}{\left( n-1\right)
			^{1-\frac{1}{\alpha }}} \\
		&=&2\frac{\left[ x_{0}\left( 1+x_{0}\right) \right] ^{\frac{1}{\alpha }}}{%
			\left( n-1\right) ^{1-\frac{1}{\alpha }}}.
	\end{eqnarray*}%
	Subcase $b).$ Assume now that $\left( j-k\right) ^{\alpha }\geq \frac{%
		k\left( n+j-1\right) }{\left( n-1\right) }$ it means $\left( j-k\right)
	^{\alpha }\left( n-1\right) -k\left( n+j-1\right) \geq 0.$ Since $n$ and $j$
	are fixed, we can define the real function 
	\begin{equation*}
		g\left( t\right) :=\left( j-t\right) ^{\alpha }\left( n-1\right) -t\left(
		n+j-1\right) ,
	\end{equation*}%
	for all $t\in \mathbb{R}$. For $t\leq \frac{n}{n+1}j$, $g\left( t\right) $
	is nonicreasing on the interval $\left[ 0,\frac{n}{n+1}j\right] .$ Really
	since%
	\begin{equation*}
		g^{\prime }\left( t\right) =-\alpha \left( j-t\right) ^{\alpha -1}\left(
		n-1\right) -\left( n+j-1\right) <0.
	\end{equation*}%
	We get%
	\begin{eqnarray*}
		g\left( \frac{nj}{n+1}\right) &=&\left( j-\frac{nj}{n+1}\right) ^{\alpha
		}\left( n-1\right) -\frac{nj}{n+1}\left( n+j-1\right) \\
		&=&\left( \frac{nj+j-nj}{n+1}\right) ^{\alpha }\left( n-1\right) -\frac{nj}{%
			n+1}\left( n+j-1\right) \\
		&=&j^{\alpha }\frac{\left( n-1\right) }{\left( n+1\right) ^{\alpha }}-\frac{n%
		}{n+1}j\left( n+j-1\right) \\
		&\leq &j^{\alpha }-j\left( n+j-1\right) =j\left( j^{\alpha -1}-\left(
		n+j-1\right) \right) <0.
	\end{eqnarray*}
	
	In the last inequality we used that $n\geq j^{\alpha -1}$ from the
	hypothesis, and $n\geq 4.$
	
	Considering the result above and because $g\left( 0\right) =j^{\alpha
	}\left( n-1\right) >0$, by the monotonicity of $g$ too, it follows that
	there exists $\tilde{k}\in \mathbb{N}$, $\tilde{k}<\frac{nj}{n+1}$ of
	maximum value, such that $g\left( \tilde{k}\right) =\left( j-\tilde{k}%
	\right) ^{\alpha }\left( n-1\right) -\tilde{k}\left( n+j-1\right) \geq 0.$
	Denoting $k_{2}=\tilde{k}+1$ and reasoning as in case $2),$ subcase $b)$ we
	have $g\left( k_{2}\right) <0.$ Furter, reasoning as in case $3)$, subcase $%
	a)$ we have $j-k_{2}<\left[ \frac{j\left( n+j-1\right) }{n-1}\right] ^{\frac{%
			1}{\alpha }}$. It follows 
	\begin{eqnarray*}
		\text{$\underline{M}$}_{\tilde{k},n,j}\left( x_{0}\right) &=&m_{\tilde{k}%
			,n,j}\left( x_{0}\right) \left( x_{0}-\frac{\tilde{k}}{n-1}\right) \leq 
		\frac{j+1}{n-1}-\frac{\tilde{k}}{n-1} \\
		&=&\frac{j-k_{2}}{n-1}+\frac{2}{n-1}<3\frac{\left[ x_{0}\left(
			1+x_{0}\right) \right] ^{\frac{1}{\alpha }}}{\left( n-1\right) ^{1-\frac{1}{%
					\alpha }}}.
	\end{eqnarray*}%
	In the light of Lemma \ref{Lemma4.4}, $(ii)$, it follows that $\underline{M}%
	_{\tilde{k},n,j}\left( x_{0}\right) \geq \underline{M}_{\tilde{k}%
		-1,n,j}\left( x_{0}\right) \geq ...\geq \underline{M}_{0,n,j}\left(
	x_{0}\right) .$ Thus we obtain 
	\begin{equation*}
		\text{$\underline{M}$}_{k,n,j}\left( x_{0}\right) <3\frac{\left[ x_{0}\left(
			1+x_{0}\right) \right] ^{\frac{1}{\alpha }}}{\left( n-1\right) ^{1-\frac{1}{%
					\alpha }}}
	\end{equation*}%
	for any $k\in \left\{ 0,1,...,\tilde{k}\right\} .$
	
	In both subcases, Lemma \ref{Lemma4.3} $(iii)$ too, we get 
	\begin{equation*}
		M_{k,n,j}\left( x_{0}\right) <6\frac{\left[ x_{0}\left( 1+x_{0}\right) %
			\right] ^{\frac{1}{\alpha }}}{\left( n-1\right) ^{1-\frac{1}{\alpha }}}.
	\end{equation*}%
	So, taking into consideration the fact that%
	\begin{equation*}
		\max \left\{ 
		\begin{array}{c}
			3\frac{\left[ x_{0}\left( 1+x_{0}\right) \right] ^{\frac{1}{\alpha }}}{%
				\left( n-1\right) ^{1-\frac{1}{\alpha }}},4\frac{\left[ x_{0}\left(
				1+x_{0}\right) \right] ^{\frac{1}{\alpha }}}{\left( n-1\right) ^{1-\frac{1}{%
						\alpha }}}, \\ 
			5\frac{\left[ x_{0}\left( 1+x_{0}\right) \right] ^{\frac{1}{\alpha }}}{%
				\left( n-1\right) ^{1-\frac{1}{\alpha }}},6\frac{\left[ x_{0}\left(
				1+x_{0}\right) \right] ^{\frac{1}{\alpha }}}{\left( n-1\right) ^{1-\frac{1}{%
						\alpha }}}%
		\end{array}%
		\right\} \leq 6\frac{\left[ x_{0}\left( 1+x_{0}\right) \right] ^{\frac{1}{%
					\alpha }}}{\left( n-1\right) ^{1-\frac{1}{\alpha }}},
	\end{equation*}%
	we have desired result.
\end{pf}

\section{Weighted Rate of Convergence}

We see that previous results works for a fixed $x_{0}$ point or finite
intervals. But if we want to obtain a uniform approximation order on
infinite intervals, then we should use weighted modulus of continuities.

Before giving useful properties about these type of modulus of continuities,
let us recall the following spaces and norm (see, for instance, \cite%
{Gadjiev1}, \cite{Gadjiev2}):%
\begin{eqnarray*}
	B_{\rho }(\mathbb{R}) &=&\left\{ f:\mathbb{R\rightarrow }\left. \mathbb{R}%
	\right\vert \text{ a constant }M_{f}\text{ depending on }f\text{ exists}%
	\right. \\
	&&\left. \text{such that }\left\vert f\right\vert \leq M_{f}\rho \right\} ,
\end{eqnarray*}%
\begin{equation*}
	C_{\rho }(\mathbb{R})=\left\{ f\in \left. B_{\rho }(\mathbb{R})\right\vert 
	\text{ }f\text{ continuous on }\mathbb{R}\right\} ,
\end{equation*}%
endowed with the norm:%
\begin{equation*}
	\left\Vert f\right\Vert _{\rho }=\sup_{x\geq 0\geq }\frac{\left\vert
		f(x)\right\vert }{\rho (x)}.
\end{equation*}

In order to obtain rate of weighted approximation of the positive linear
operators defined on infinite intervals, various weighted modulus of
continuities are introduced. Some of them include term $h$ in the
denominator of the supremum expression. In the chronological order, let us
refer to some related papers as \cite{Achieser}, \cite{Freud}, \cite{Amanov}%
, \cite{Gadjieva}, \cite{Dogru7}, \cite{Moreno}, \cite{Aral}, \cite{Holhos1}.

The weighted modulus defined in \cite{Achieser}, in order to obtain weighted
approximation properties of Szasz-Mirakyan operators on $\mathbb{R}_{+}$.

In \cite{Gadjieva}, second author, jointly with Gadjieva was introduced the
following modulus of continuity:%
\begin{equation}
	\Omega (f;\delta )=\sup_{0\leq x,\left\vert h\right\vert \leq \delta }\frac{%
		\left\vert f(x+h)-f(x)\right\vert }{(1+h^{2})(1+x^{2})}.  \label{5.1}
\end{equation}%
There are some papers including rates of weighted approximation with the
help of $\Omega (f;\delta ).$ (see, for instance, \cite{Karsli}, \cite{Deniz}%
, \cite{Dogru9} and \cite{Ispir}).

And then in \cite{Dogru7}, second author defined the following modulus of
continuity:%
\begin{equation}
	\omega _{\rho }(f;\delta )=\sup_{0\leq x,\left\vert h\right\vert \leq \delta
	}\frac{\left\vert f(x+h)-f(x)\right\vert }{\rho (x+h)}  \label{5.2}
\end{equation}%
where $\rho (x)\geq \max (1,x).$

In \cite{Dogru7}, the author was introduced a generalization of the
Gadjiev-Ibragimov operators which includes many well-known operators and
obtain its rate of weighted convergence with the help of $\omega _{\rho
}(f;\delta )$ defined in (\ref{5.2}).

In \cite{Moreno}, Moreno introduced another type of modulus of continuity in
(\ref{5.2}) as follows%
\begin{equation*}
	\overline{\Omega }_{\alpha }(f;\delta )=\sup_{0\leq x,\left\vert
		h\right\vert \leq \delta }\frac{\left\vert f(x+h)-f(x)\right\vert }{%
		1+(x+h)^{\alpha }}.
\end{equation*}%
In \cite{Aral}, Gadjiev and Aral defined the following modulus of continuity:%
\begin{equation*}
	\widetilde{\Omega }_{\rho }(f;\delta )=\sup_{x,t\in \mathbb{R}%
		_{+},\left\vert \rho (t)-\rho (x)\right\vert \leq \delta }\frac{\left\vert
		f(t)-f(x)\right\vert }{\left( \left\vert \rho (t)-\rho (x)\right\vert
		+1\right) \rho (x)}
\end{equation*}%
where $\rho (0)=1$ and $\inf_{x\geq 0}\rho (x)\geq 1.$

It is obvious that by choosing $\alpha =2,$ in the definition of $\overline{%
	\Omega }_{\alpha }(f;\delta ),$ then we obtain $\overline{\Omega }%
_{2}(f;\delta )=\omega _{\rho _{0}}(f;\delta )$ for $\rho _{0}(x)=1+x^{2},$
and if we choose $\alpha =2+\lambda $ in the definition of $\overline{\Omega 
}_{\alpha }(f;\delta ),$ then we obtain%
\begin{equation*}
	\widehat{\Omega }_{\rho \lambda }(f;\delta )=\sup_{0\leq x,\left\vert
		h\right\vert \leq \delta }\frac{\left\vert f(x+h)-f(x)\right\vert }{%
		1+(x+h)^{2+\lambda }}
\end{equation*}%
(see \cite{Agratini}).

Finally, in \cite{Holhos1}, Holho\c{s} defined a more general weighted
modulus of continuity as%
\begin{equation*}
	\omega _{\varphi }(f;\delta )=\sup_{0\leq x\leq y,\left\vert \varphi
		(y)-\varphi (x)\right\vert \leq \delta }\frac{\left\vert
		f(x)-f(y)\right\vert }{\rho (x)+\rho (y)}
\end{equation*}%
such that, for $\varphi (x)=x,$ this modulus of continuity is equivalent to $%
\Omega (f;\delta )$ defined in (\ref{5.1})$.$

Also, let $C_{\rho }^{0}(\mathbb{R})$ be the subspace of all functions in $%
C_{\rho }(\mathbb{R})$ such that $\lim_{\left\vert x\right\vert \rightarrow
	\infty }\frac{f(x)}{\rho (x)}$ exists finitely.

Notice also that some remarkable properties about these type of modulus of
continuities can be found in \cite{DL}.

In the light of these definitions, we can give the following theorem:

\begin{thm}
	\label{Theorem6.1} Let $f:\left[ 0,\infty \right) \rightarrow \mathbb{R}_{+}$
	be continuous. Then for all $x\in $ $[0,\infty ),$ which also satisfy $%
	x_{0}^{\alpha -2}\leq n-1$ and $k<nx$ and $n\geq j^{\alpha -1},$ $n\geq 4,$
	we have the following order of approximation for the operators (\ref{2.1})
	to the function $f$ by means of the weighted modulus of continuity defined
	in (\ref{5.2}). Then for each $f\in C_{\rho _{0}}^{0}(\mathbb{R}_{+})$, we
	have%
	\begin{equation}
		\frac{\left\vert V_{n}^{(M)}(f)(x)-f(x)\right\vert }{(\rho _{0}(x))^{2}}\leq %
		\left[ 
		\begin{array}{c}
			\frac{(1+9x^{2})\left( 1+6\left[ x\left( 1+x\right) \right] ^{\frac{1}{%
						\alpha }}\right) }{(1+x^{2})^{2}} \\ 
			\text{ }\omega _{\rho _{0}}\left( f;\frac{1}{\left( n-1\right) ^{1-\frac{1}{%
						\alpha }}}\right) 
		\end{array}%
		\right] ,\text{ }  \label{5.3}
	\end{equation}%
	for all $n\in \mathbb{N},$ $n\geq 4$ where $\rho _{0}(x)=1+x^{2}$ and $%
	\alpha =2,3,...$ .
\end{thm}

\begin{pf}
	By using the properties of $\omega _{\rho _{0}}(f;\delta ),$(see \cite%
	{Moreno} ), we can write 
	\begin{equation}
		\left\vert V_{n}^{(M)}(f)(x)-f(x)\right\vert \leq \left[ 
		\begin{array}{c}
			\left( 1+(2x+V_{n}^{(M)}\left( e_{1}\right) \left( x\right) )^{2}\right)  \\ 
			\left( \frac{1}{\delta }V_{n}^{(M)}\left( \varphi _{x}\right) \left(
			x\right) +1\right) \omega _{\rho _{0}}(f;\delta )%
		\end{array}%
		\right] .  \label{5.4}
	\end{equation}%
	In the proof of Theorem \ref{Theorem5.1}, we obtain%
	\begin{equation}
		V_{n}^{(M)}\left( \varphi _{x}\right) \left( x\right) \leq 6\frac{\left[
			x\left( 1+x\right) \right] ^{\frac{1}{\alpha }}}{\left( n-1\right) ^{1-\frac{%
					1}{\alpha }}},  \label{5.5}
	\end{equation}
	for all $n\in \mathbb{N}$, $n\geq 4$ and $x\in \left[ 0,\infty \right) ,$ $%
	x_{0}^{\alpha -2}\leq n-1.$ On the other hand, we have%
	\begin{eqnarray*}
		V_{n}^{(M)}\left( e_{1}\right) \left( x\right)  &=&\frac{\bigvee%
			\limits_{k=0}^{\infty }\frac{\left( n+k-1\right) !}{k!\left( n-1\right) !}%
			\frac{x^{k}}{\left( 1+x\right) ^{n+k}}\frac{k}{n}}{\bigvee\limits_{k=0}^{%
				\infty }\binom{n+k-1}{k}\frac{x^{k}}{\left( 1+x\right) ^{n+k}}}=\frac{%
			x\bigvee\limits_{k=1}^{\infty }\binom{n+k-1}{k-1}\frac{x^{k-1}}{\left(
				1+x\right) ^{n+k}}}{\bigvee\limits_{k=0}^{\infty }\binom{n+k-1}{k}\frac{x^{k}%
			}{\left( 1+x\right) ^{n+k}}} \\
		&=&\frac{x\bigvee\limits_{k=0}^{\infty }\binom{n+k}{k}\frac{x^{k}}{\left(
				1+x\right) ^{n+k+1}}}{\bigvee\limits_{k=0}^{\infty }\binom{n+k-1}{k}\frac{%
				x^{k}}{\left( 1+x\right) ^{n+k}}}.
	\end{eqnarray*}%
	And also, since $k<nx$ then we obtain%
	\begin{equation}
		V_{n}^{(M)}\left( e_{1}\right) \left( x\right) \leq x.  \label{5.6}
	\end{equation}%
	So, using the inequalities (\ref{5.5}) and (\ref{5.6}) in (\ref{5.4}) and
	choosing%
	\begin{equation*}
		\delta =\frac{1}{\left( n-1\right) ^{1-\frac{1}{\alpha }}},
	\end{equation*}%
	the proof is completed.
\end{pf}

This theorem allows us to express the following weighted uniform
approximation result.

\begin{thm}
	\label{Theorem6.2} Let $f:\left[ 0,\infty \right) \rightarrow \mathbb{R}_{+}$
	be continuous. Then for all $x\in $ $[0,\infty ),$ which also satisfy $%
	x_{0}^{\alpha -2}\leq n-1$ and $k<nx$ and $n\geq j^{\alpha -1},$ $n\geq 4,$
	we have the following order of approximation for the operators (\ref{2.1})
	to the function $f$ by means of the weighted modulus of continuity defined
	in (\ref{5.2}). Then for each $f\in C_{\rho _{0}}^{0}(\mathbb{R}_{+})$, we
	have%
	\begin{equation}
		\left\Vert V_{n}^{(M)}(f)(x)-f(x)\right\Vert _{\rho _{0}^{2}(x)}\leq 70\text{
		}\omega _{\rho _{0}}\left( f;\frac{1}{\left( n-1\right) ^{1-\frac{1}{\alpha }%
		}}\right) ,\text{ }  \label{5.7}
	\end{equation}%
	for all $n\in \mathbb{N},$ $n\geq 4$ where $\rho _{0}(x)=1+x^{2}$ and $%
	\alpha =2,3,...$ .
\end{thm}

\begin{pf}
	By using the inequalities $\frac{1}{1+x^{2}}\leq 1,$ $\frac{x^{2}}{1+x^{2}}%
	\leq 1,$ and $\frac{\left[ x\left( x+1\right) \right] ^{\frac{1}{\alpha }}}{%
		1+x^{2}}\leq 1$\ we have%
	\begin{equation}
		\frac{(1+9x^{2})\left( 1+6\left[ x\left( 1+x\right) \right] ^{\frac{1}{%
					\alpha }}\right) }{(1+x^{2})^{2}}\leq 70.  \label{5.8}
	\end{equation}%
	If we use (\ref{5.8}) in (\ref{5.3}), we obtain desired result.
\end{pf}

\begin{rem}
	So, Theorem \ref{Theorem5.1}, Theorem \ref{Theorem6.1} and Theorem \ref%
	{Theorem6.2} show that the orders of pointwise approximation, weighted
	approximation and weighted uniform approximation are $1/\left( n-1\right)
	^{1-\frac{1}{\alpha }}.$ For big enough $\alpha $ , $1/\left( n-1\right) ^{1-%
		\frac{1}{\alpha }}$ tends to $1/(n-1)$. As a result, since $1-\frac{1}{%
		\alpha }\geq \frac{1}{2}$ for $\alpha =2,3,...$ this selection of $\alpha $
	improving the order of approximation.
\end{rem}

\begin{ack}
	{Dedicated to Professor Abdullah Alt\i n on the occasion of his
		75 th birthday, with high esteem.}
\end{ack}

\end{document}